\documentclass[11pt,twoside]{article}
\usepackage{amsfonts,amsmath,epsfig,amssymb,amsthm}

\usepackage{cite}

\usepackage{xcolor}

\numberwithin{equation}{section}
\newfont{\ctv}{msam10}
\newcommand{\bbox}{\mbox{\ctv \symbol{3}}}
\def\QED{{$\hfill\bbox$}}
\newenvironment{pf}[1]{\par\vskip1mm{\noindent\it #1.}\ }{\QED\par\vskip2mm}
\newtheorem{theorem}{Theorem}[section]

\newtheorem{lemma}{Lemma}[section]
\newtheorem{definition}{Definition}[section]

\newtheorem*{theoremH}{Hypotheses (H)}
\newtheorem*{theoremU}{Hypotheses (U)}

\newcommand{\Barint}{-\kern-.15in\int}
\newcommand{\barint}{-\kern-.12in\int}
\renewcommand{\r}{\mathbb{R}}
 
\DeclareMathOperator{\haus}{haus} 
\DeclareMathOperator{\co}{co} 

\def\bpf{\begin{pf}}
\def\epf{\end{pf}}

\begin{document}

\thispagestyle{empty}
 \setcounter{page}{1}

 \title{Relaxation in population dynamics models with hysteresis
 \footnote{The research of the first author was supported by the Program for Innovative Research Team in Science and Technology in Fujian Province University, by Quanzhou High-Level Talents Support Plan under Grant 2017ZT012,
by the Scientific Research Fund of Huaqiao University no. 605-50Y19017, and by RFBR grant no. 18-01-00026}
}
 \author{Sergey A. Timoshin\footnote{Fujian Province University Key Laboratory of Computational Science, School of Mathematical Sciences, Huaqiao University, Quanzhou 362021, China -- and --
Matrosov Institute for System Dynamics and Control
 Theory,
 Russian Academy of Sciences, Lermontov str. 134,
664033  Irkutsk,
 Russia, E-mail: {\tt
sergey.timoshin@gmail.com}. } $^{,}$\footnote{Corresponding
author.}
\ and Toyohiko Aiki\footnote{Department of Mathematical and
Physical Sciences, Faculty of Science, Japan Women's University,
2-8-1 Mejirodai, Bunkyo-ku, Tokyo, 112-8681, Japan, E-mail: {\tt
aikit@fc.jwu.ac.jp}.} }

\date{\emph{submitted on} August 6, 2019}

\maketitle \hyphenation{Ca-ra-the-odo-ry}

\begin{abstract}
\noindent {The present paper is concerned with a nonlinear partial differential control system subject to a state-dependent and nonconvex control constraint. This system models the dynamics of populations in the vegetation--prey--predator framework and takes account of diffusive and hysteresis effects appearing in the process. We prove the existence of solutions to our system and show that they are close in a suitable sense to solutions of the system with the convexified control constraint.}
\end{abstract}

\small\textbf{Keywords}:   biological diffusion models, hysteresis, evolution control systems, state-dependent constraints, relaxation.


\vspace{0.3cm}

\small\textbf{2010 Mathematics Subject Classification}:  49J45,
49J53, 93C20, 93C73.

\section{Introduction} \label{number of section:introduction}

The motivation of the present paper comes from problems arising in  prey--predator models when diffusive effects in the dynamics of the prey and predator populations are taken into account and
the evolution of the food density for the prey exhibits a hysteretic character. Aiming at achieving a possible optimization of the population dynamics process by way of controlling the growth rate of the prey, we introduce the following dynamical control problem:

\begin{equation}
\sigma_t -av_t + \partial I_{v,w}(\sigma) \ni F(\sigma,v,w)  \qquad \text{in} \; Q(T)
,\label{1.1}
\end{equation}
\begin{equation} v_t - \Delta v=h(\sigma, v,w)\,u \qquad \text{in}
\; Q(T) , \label{1.2}\end{equation}
\begin{equation} w_t - \Delta w=g(\sigma, v,w)  \qquad \text{in}
\; Q(T) , \label{1.3}\end{equation}
\begin{equation}
\sigma(x,0)=\sigma_0(x), \quad v(x,0)=v_0(x), \quad w(x,0)=w_0(x)   \quad \text{in} \; \Omega ,\label{1.4}
\end{equation}
\begin{equation}
\;\frac{\partial v}{\partial n}=\frac{\partial w}{\partial n}=0 \qquad \text{on} \; \partial\Omega\times [0,T] .\label{1.5}
\end{equation}

Here,  $Q(T) := [0,T]\times \Omega$ with $T>0$ being a fixed final time and $\Omega$ a bounded domain in $\mathbb{R}^N$, $N\leq 3$,  with smooth
boundary $\partial\Omega$. $I_{v,w}(\cdot)$ is the
indicator function of the interval $[f_*(v,w),
f^*(v,w)]$, $\partial I_{v,w}(\cdot)$ is
its subdifferential in the sense of convex analysis.  $f_*,f^*:\mathbb{R}^2\to \mathbb{R}$, $F, h, g:\mathbb{R}^3\to \mathbb{R}$ are prescribed functions
with properties enlisted in the next section, $a$ is a constant, $\sigma_0$, $v_0$,  $w_0$ are given  initial
conditions, and $\partial/\partial n$ is the outward normal derivative on $\partial\Omega$. The function $u$ on the right-hand side of $(\ref{1.2})$ plays the role of control.

\hspace{0.6cm} In our model, the unknown variables $\sigma$,  $v$, and $w$ represent the densities of
the food for the prey (vegetation), the prey
and the predator, respectively. The evolution of the food density is characterized by a hysteretic relationship with the hysteresis region generated by the characteristic curves $\sigma=f_*(v,w)$ and $\sigma=f^*(v,w)$ (cf. \cite{Visintin})
describing the situation in which the growth rate of the food for the prey depends not only on the present state of preys and predators, but also on their immediately preceding density history. This hysteretic behaviour is captured by introducing a hysteresis operator into the system. In its turn, the latter is represented by adding  the subdifferential term to Eq. $(\ref{1.1})$.   At this point, we would like also to mention that every scalar
return point memory hysteresis operator can be represented by first order differential inclusion with a one-parameter family of
indicator functions (see \cite[Theorem 2.7.7]{BrokateSprekels}).

\hspace{0.6cm} Nonlinear phenomena of hysteresis type are encountered in many branches of natural and applied sciences ranging from the physics of materials to economics. The biological literature has also repeatedly described the situation when the way the state variables of a process change after the system's parameters have been changed is different from the way the variables change back when the parameters regain their former values and a hysteresis loop is thus formed. Note, however, that, to the authors' knowledge, contributions with the rigorous mathematical treatment of biological processes with hysteresis, let alone controlled biological systems with hysteresis, are still very few in number (see, e.g., \cite{VisintinSIAM, AikiKopfova, Gurevich, AikiMinchev, scl2019}).

\hspace{0.6cm} The famous spruce budworm population dynamics models  can serve an example of practical situation where our results may find potential applications. These models describe the budworm--forest ecosystem consisting of a forest insect pest (spruce budworm) endemic to eastern North America which defoliates balsam fir and several other tree species in the boreal forest and is a prey for avian predators of the forest. The modelling and the subsequent study of the budworm--forest interactions are very important in the forest ecology as the budworm during its outbreaks causes substantial damage to the forest destroying a large number of trees (see \cite{Murray, Ludwig} for a particular instance of spruce budworm dynamics modelling).

\hspace{0.6cm} System $(\ref{1.1})$--$(\ref{1.5})$ is considered subject to the following state-dependent control constraint
\begin{equation}
u\in U(t,x,\sigma,v,w) \qquad \text{in} \; Q(T) ,\label{1.6}
\end{equation}

where the control constraint multifunction $U:Q(T)\times \mathbb{R}^3\to \mathbb{R}$ has compact, but not necessarily convex values. We note that while the nonconvexity of values of $U$ might be a biologically relevant assumption, it poses certain difficulties for mathematical and numerical analysis of the control problem.
Hence, along with $(\ref{1.6})$  we consider the following alternative (convexified) control constraint
\begin{equation}
u\in \co U(t,x,\sigma,v,w) \qquad \text{in} \; Q(T) ,\label{1.7}
\end{equation}

where $\co U$ denotes the convex hull of the set
$U$, which is the smallest under inclusion convex set containing $U$.  The corresponding systems $(\ref{1.1})$--$(\ref{1.6})$ and $(\ref{1.1})$--$(\ref{1.5})$, $(\ref{1.7})$ will in the sequel be
refereed to and denoted as the given (or original) $(P)$ and convexified (or
relaxed) $(RP)$ problems, respectively.

\hspace{0.6cm} The main aim of the present paper then is to establish the \emph{existence} of solutions to Problems $(P)$ and $(RP)$ and to show that the solutions of the two problems are close in a prescribed sense. Namely, we establish the so-called \emph{relaxation} property for system $(P)$ asserting that its solutions are dense in an appropriate topology among the solutions of system $(RP)$. The exact meaning in which solutions to Problems $(P)$ and $(RP)$ and the relaxation property are understood is explained in the next section.

\hspace{0.6cm} We note that control system $(\ref{1.1})$--$(\ref{1.6})$  is a modification of the following control system, coupled with the relevant initial boundary conditions and  control constraint, considered recently in \cite{scl2019} to describe the evolution of populations in the prey-predator framework when diffusion of the vegetation is being accounted for:
\begin{equation}
\sigma_t -(\lambda(v))_t -\kappa \Delta\sigma+ \partial I_{v,w}(\sigma) \ni F(\sigma,v,w)\,u   \qquad \text{in} \; Q(T)
,\label{1.8}
\end{equation}
\begin{equation} v_t - \Delta v=h(\sigma, v,w)\qquad \text{in}
\; Q(T) , \label{1.9}\end{equation}
\begin{equation} w_t - \Delta w=g(\sigma, v,w)  \qquad \text{in}
\; Q(T) , \label{1.10}\end{equation}

where $\lambda:\mathbb{R}\to\mathbb{R}$ is a given function and $\kappa>0$ is a diffusion parameter. In \cite{scl2019} we proved the existence of solutions for this control problem.

\hspace{0.6cm}  There are a number of reasons  for considering our control problem in the form $(\ref{1.1})$--$(\ref{1.3})$ in place of $(\ref{1.8})$--$(\ref{1.10})$. First, a seemingly simplifying assumption that $\kappa=0$ renders, in actual fact, the mathematical investigation of system $(\ref{1.1})$--$(\ref{1.3})$ more challenging as, in this case, less spatial regularity of the state $\sigma$ is entailed and the dependence of $\sigma$ on $x$ may not be necessarily smooth. On the other hand, the absence of vegetation diffusion is quite natural in many biological models, in particular, in the spruce budworm population dynamics model mentioned above when considered on the short-to-mid term timescale. Second, the inclusion of  the external controller $u$ to the second equation of the system as in $(\ref{1.1})$--$(\ref{1.3})$ instead of the first one as in $(\ref{1.8})$--$(\ref{1.10})$ is more justifiable from an ecological viewpoint as the typical controlling actions available usually directly affect the rate of change in budworm population, e.g. direct spraying of insecticides, removal of infected trees and so on. The price we need to pay for the above ameliorations to the model is that we are constrained to consider only the case of a linear function  $\lambda$. Note, however, that this is not a real restriction from the biological perspective as the function $\lambda$ considered in typical examples from the population dynamics is linear. Moreover, when $w$ is fixed and $F\equiv 0$, $a=1$ Eq. $(\ref{1.1})$ recovers the differential representation of the generalized stop operator (cf. \cite{Visintin}).

\hspace{0.6cm}   In conclusion, we mention that when considering optimal control problems, necessary optimality conditions are usually obtained only for convex problems (convex cost functional and convex constraints). At the same time, numerical algorithms for optimal control problems are largely based on necessary optimality conditions. In this respect, our relaxation results provide a step towards justification of the passage from real life nonconvex problems to amenable to calculations convex problems.

\hspace{0.6cm} At the end of the introduction, we mention that some optimal control problems with nonconvex control constraints have been recently considered in connection with fractional calculus \cite{Debbouche2017, Debbouche2018} and stochastic analysis \cite{Debbouche2020}. In this respect, a combination of fractional and/or stochastic calculus with hysteresis systems might prove to yield models better reflecting the properties of real-life problems thus opening a new perspective direction of research.

\section{Notation and assumptions} \label{section2}

Denote by $H$ the Hilbert space $L^2(\Omega)$ with
the usual scalar product $\langle\cdot,\cdot\rangle_H$ and the norm $|\cdot|_H$, and let $V$
 be the Sobolev space $H^1(\Omega$) equipped with the norm
$|v|_V=\langle v,v\rangle_V^{1/2},$ where $\langle
v,w\rangle_V=\langle v,w\rangle_H+\int_{\Omega} \langle \nabla v(x), \nabla w(x)
\rangle_{{\mathbb R}^N} \,dx$, $v,w\in V$. Let $V'$ be the dual space of $V$ and $\langle \cdot ,\cdot
\rangle$ stand for the duality pairing between $V'$ and $V$.
Define the operator $\label{} -\Delta_N: \; D(-\Delta_N)\subset H\to H$
as the restriction of the linear continuous operator
${\mathcal R}:V\to V'$, $\langle {\mathcal R}v,w\rangle =\int_{\Omega} \langle \nabla v(x), \nabla w(x)
\rangle_{{\mathbb R}^N} \,dx$, $v,w\in V$,
 to the subset of $V$ consisting
of the elements $v$ such that ${\mathcal R}v \in H$. Then, we have
$$
D(-\Delta_N)= \left\{v\in H^2(\Omega); \; {\partial v}/{\partial n}
= 0 \; \; \hbox{in} \;\; H^{1/2}(\partial \Omega ) \right\}
$$
and
\begin{equation*}
\label{} -\Delta_N v=-\Delta v \quad \hbox{for all} \;\; v\in
D(-\Delta_N).
\end{equation*}

\hspace{0.6cm} Given a convex, lower semicontinuous
function
$\varphi:H\rightarrow\mathbb{R\cup\{+\infty\}}$ which is not identically $+\infty$, its subdifferential $\partial
\varphi(x)$ at a point $x\in H$ is the set
$$\partial
\varphi(x)=\{h\in H : \langle h, y-x\rangle_H\leq
\varphi(y)-\varphi(x) , \; \forall y\in H\} .$$ The subdifferential mapping
$\partial \varphi:H\to H$ is a maximal monotone operator. A
multivalued operator $A:H\rightarrow H$ is said to be monotone if for
any $x_i \in \mbox{dom\,} A:=\{x\in H:
Ax\neq\emptyset\}$, and any $h_i\in Ax_i$, $i=1,2$, the
inequality $\langle x_1-x_2, h_1 -h_2\rangle\geq 0$ holds.

\hspace{0.6cm} For a Banach space $X$ we denote by $d_X(x,C)$ the distance from a point $x\in X$ to a set
$C\subset X$. Then, the Hausdorff metric  on the space of closed
bounded subsets of $X$, denoted $cb (X)$, is the function:
$$\haus_X(C,D)=\max\{\sup\limits_{x\in C}d_X(x,D),\sup\limits_{y\in D}d_X(y,C)\},
\quad C, D \in cb(X).$$

\hspace{0.6cm} Given a metric space $Y$ and a point $y_0\in Y$, a multivalued mapping $A:Y\to X$ is called lower semicontinuous at $y_0$ if for any $x_0\in A(y_0)$ and any sequence $y_n\in Y$, $n\geq 1$, converging to $y_0$, there exists a sequence $x_n\in A(y_n)$, $n\geq 1$, converging to $x_0$. The mapping $A$ is lower semicontinuous on a subset of $Y$ if it is lower semicontinuous at every point of this subset.

\hspace{0.6cm} A multivalued mapping $A$ from a measurable
space $({\mathcal E},{\mathcal A})$ to
$cb(X)$ is called measurable if $\{\tau \in {\mathcal E}; \;
A(\tau )\cap C\neq \emptyset \}\in {\mathcal A}$ for any closed
set $C\subset X$.

\hspace{0.6cm} We introduce now the hypotheses on the data of our Problem $(P)$. These hypotheses are valid throughout the rest of the paper.

\begin{theoremH} \*
\begin{itemize}
\item[\textup{\textbf{(H1)}}] the functions $f_*, f^*\in C^2(\mathbb{R}^2)\cap
W^{2,\infty}(\mathbb{R}^2)$ are such that $0\leq f_*\leq f^*\leq 1$ on
$\mathbb{R}^2$;

\item[\textup{\textbf{(H2)}}] the functions $F, h, g:\mathbb{R}^3\to \mathbb{R}$ are
Lipschitz continuous (with a common Lipschitz constant $L>1$) and are such that $h(\sigma,0,w)=0$ for $\sigma\in [0,1]$, $w\in \mathbb{R}$, $g(\sigma,v,0)=0$ for $\sigma\in [0,1]$, $v\in \mathbb{R}$;

\item[\textup{\textbf{(H3)}}] the initial conditions $\sigma_0, v_0, w_0\in L^{\infty}(\Omega)\cap V$ are such that $v_0\geq
0$, $w_0\geq 0$ and $f_*(v_0,w_0)\leq \sigma_0\leq f^*(v_0,w_0)$ a.e. on $\Omega$.

\end{itemize}
\end{theoremH}

\hspace{0.6cm} With respect to the bounds in the first hypothesis above we note that the fact that the vegetation $\sigma$
is constant ($=1$ after rescaling) when the prey population $v$ is zero, and $\sigma=0$
 if $v$ exceeds a certain critical value is a
natural assumption from a biological viewpoint (see also Definition 2.1 $(iii)(a)$ below).

\hspace{0.6cm} The next hypothesis lists the assumptions we impose on the control constraint $(\ref{1.6})$.

\begin{theoremU} The multivalued mapping $U:[0,T]\times \Omega \times {{\mathbb R}}\times {{\mathbb R}}\times {{\mathbb R}}
\to {cb}({\mathbb R})$ has the following properties:

\begin{itemize}
\item[\textup{\textbf{(U1)}}] the mapping $(t,x)\to U(t,x,\sigma,v,w), \; \sigma,v,w\in {\mathbb R},$
is measurable;

\item[\textup{\textbf{(U2)}}] there exists a constant $m>0$ such that
\begin{equation*} \label{}
  |U(t,x,\sigma,v,w)| \leq m \quad
\mbox{a.e. on} \;\, Q(T), \sigma,v,w \in \mathbb{R};
\end{equation*}

\item[\textup{\textbf{(U3)}}]  there exists  $k\in
L^2(0,T;{{\mathbb R}}^+)$ such that
\begin{align*} \label{}
 \haus_{\mathbb R}(U(t,x,\sigma_1,v_1,w_1)&,U(t,x,\sigma_2,v_2,w_2)) \\ &\leq k(t) (|\sigma_1-\sigma_2|+|v_1-v_2|+|w_1-w_2|)
\end{align*}

a.e. on $Q(T)$, $\sigma_i,v_i,w_i\in \r$, $i=1,2$.

\end{itemize}
\end{theoremU}

\hspace{0.6cm} In order to define a solution to our problems $(P)$ and $(RP)$ we first define the multivalued mapping
\begin{equation*}
\label{} {\mathcal U}(t,\sigma,v,w)=\{u\in H; \; u(x)\in U(t,x,\sigma(x),v(x),w(x))
\; \; \hbox{a.e. on} \; \; \Omega \}, \quad  \sigma, v, w \in H,
\end{equation*}

and the set
$${\mathcal{K}}(v,w)=\{\sigma\in H; \; f_*(v(x),w(x))\leq \sigma(x) \leq f^*(v(x),w(x)) \;\;\;\text{a.e. on}\; \Omega\} , \quad v,w\in H.$$

Then, from \cite[Lemma 3.1]{TCON2018}) we see that the following
properties hold for the mapping
${\mathcal{U}}:[0,T]\times H \times H \times H\to {cb}(H)$:
\begin{itemize}
\item[$\textbf{(}\mathcal{U}\textbf{1)}$] the mapping $t\mapsto
    {\mathcal{U}}(t,\sigma,v,w)$ is measurable, $\sigma,v,w\in
    H$;

\item[$\textbf{(}\mathcal{U}\textbf{2)}$] $|{\mathcal{U}}(t,\sigma,v,w)|_{H}\leq
    m$\;\; a.e. on $[0,T]$,  $\sigma,v,w\in
    H$, where $m>0$ is as above;

\item[$\textbf{(}\mathcal{U}\textbf{3)}$] \*
\vspace{-0.7cm}\begin{align*}
\hspace{-1.3cm}    \haus_H({\mathcal{U}}(t,\sigma_1,v_1,w_1)&,{\mathcal{U}}(t,\sigma_2,v_2,w_2))  \\ &\leq
    k(t)(|\sigma_1-\sigma_2|_H+|v_1-v_2|_H+|w_1-w_2|_H) \end{align*}
    a.e. on $[0,T]$,  $\sigma_i,v_i,w_i\in H$, $i=1,2$ for $k\in L^2(0,T;\mathbb{R}^+)$ as above.
\end{itemize}

\begin{definition} A quadruple  $\{\sigma,v,w,u\}$ is called a solution of
control system $(P)$  if
\begin{itemize}
\item[$(i)$] $\sigma\in W^{1,2}(0,T;H)$, $v,w\in W^{1,2}(0,T;H)\cap
L^\infty(0,T;V)\cap L^2(0,T;H^2 (\Omega))$;

\item[$(ii)$] $u\in L^2(0,T;H)$;

\item[$(iii)$] $ \sigma'-a v'+ \partial I_{\mathcal{K}(v,w)}(\sigma)
    \ni F(\sigma,v,w)$ \; in $H$ a.e. on $[0,T]$;

\item[$(iv)$] $ v'- \Delta_N v=h(\sigma,v,w)\,u$
\; in $H$ a.e. on $[0,T]$;

\item[$(v)$] $ w'-\Delta_N w=g(\sigma,v,w)$
\; in $H$ a.e. on $[0,T]$;

\item[$(vi)$] $ \sigma(0) = \sigma_0, \; v(0) = v_0, \; w(0) = w_0$ \; in  $H$;

\item[$(vii)$] $u(t)\in \mathcal{U}(t,\sigma(t),v(t),w(t))$ \; in $H$
for a.e. $t\in [0,T]$,

\end{itemize}

where the prime denotes the derivative with respect to $t$.

\hspace{0.6cm} A solution of control system $(RP)$ is defined
similarly replacing the last inclusion with
\begin{equation*}
u(t)\in \overline{\co} \;\mathcal{U}(t,\sigma(t),v(t),w(t)) \quad \text{in} \;
H \;\; \text{for a.e.} \; t\in [0,T] .\label{}
\end{equation*}
\end{definition}

\hspace{0.6cm} When $u$ is fixed in some appropriate set the notion of a solution for system $(\ref{1.1})$--$(\ref{1.5})$ naturally extends from Definition 2.1. So, in this case, a solution is a triple $\{\sigma,v,w\}$
satisfying $(i)$, $(iii)$--$(vi)$ of Definition 2.1 (see Theorem 3.1 of the next section).

\hspace{0.6cm} Remark that inclusion $(iii)$ in Definition 2.1 implies the following:
\begin{itemize}

\item[$(iii)(a)$] $f_*(v,w)\leq \sigma \leq f^*(v,w)$ a.e.
    in $Q(T)$;

\item[$(iii)(b)$] $ (\sigma'(t)-a v'(t)- F(\sigma(t),v(t),w(t)),
    \sigma(t)-z)_H\leq 0\,$ for all $z\in H$ with
    $f_*(v(t),w(t))\leq z \leq f^*(v(t),w(t))$ a.e. in
    $\Omega$ for a.e. $t\in [0,T]$.
\end{itemize}

\hspace{0.6cm} Given Hypotheses $(H)$ and $(U)$, the main purpose of
this work is to prove the following result.

\begin{theorem} Control systems
$(P)$ and $(RP)$ have solutions. Moreover,  for any solution
$\{\sigma_*,v_*,w_*,u_*\}$ of the latter system there exists a sequence of solutions
$\{\sigma_k,v_k,w_k,u_k\}$, $k\geq 1$, of the former one such that
$\{\sigma_k,v_k,w_k\}\to \{\sigma_*,v_*,w_*\}$ in $C([0,T];H\times H\times H)$ and $u_k\to u$ weakly
in  $L^2(0,T;H)$.
\end{theorem}

\hspace{0.6cm} We note that this last property is commonly refereed
to as \emph{relaxation}.

\section{Control-to-state solution operator} \label{section3}

The bound from Hypothesis $(U2)$ for the controls of Problem $(P)$ obviously extends to those of the convexified problem $(RP)$. In particular, all the controls of both problems belong to the set
\begin{equation}S_m=\{u\in L^2(0,T;H); \;  |u(t,x)|  \leq m \quad \mbox{a.e. on} \; Q(T)\}.
\label{3.1}
\end{equation}

We have the following theorem.

\begin{theorem} \label{th3-1}
 For any fixed $u\in S_m$ system
$(\ref{1.1})$--$(\ref{1.5})$ has a unique solution. Moreover, for
any solution $\{\sigma,v,w\}$ of $(\ref{1.1})$--$(\ref{1.5})$ with  $u\in S_m$ the following a priori estimates uniform with respect to $u$ hold
\begin{equation}\hspace{-0.5cm}0 \leq \sigma,v, w
\leq R_0 \quad \mbox{ a.e. on } Q(T), \label{3.2}
\end{equation}
\begin{align}\nonumber |\sigma' |_{L^2(0,T;
H)}&+ | v' |_{L^2(0,T; H)} + |w' |_{L^2(0,T;
H)}\\ & + |\Delta v|_{L^2(0,T; H)}+ |\Delta w|_{L^2(0,T; H)}  \label{3.3} \\ & \nonumber + | \nabla v|_{L^{\infty}(0,T; H)}+ | \nabla w|_{L^{\infty}(0,T; H)}
\leq R_0
\end{align}
for a  constant $R_0$ independent of $u$.
\end{theorem}

\bpf{Proof} The existence of a  unique solution to
$(\ref{1.1})$--$(\ref{1.5})$ for a fixed $u\in S_m$ as well as the estimate (\ref{3.2}) follow from \cite[Theorems 3.1, 3.2, and 3.10]{AikiMinchev}.

\hspace{0.6cm} We note that the bound $(\ref{3.2})$ allows us to assume that the functions $F, g, h$ are bounded on $\mathbb{R}^3$. Indeed, it is enough to restrict our analysis to the set $\{0\leq \sigma, v, w\leq R_0\}$.

\hspace{0.6cm} To derive the energy estimates (\ref{3.3}), first we multiply Eq. $(iv)$ in Definition 2.1 by $v'$ and Eq. $(v)$ in Definition 2.1 by $w'$, add the resulting equalities and invoke Young's inequality to obtain
\begin{equation} |v'|^2_H +|w'|^2_H +\frac{d}{dt}|\nabla v|^2_H+\frac{d}{dt}|\nabla w|^2_H\leq C_1   \label{3.4}\end{equation}

a.e. on $(0,T)$, where $C_1=m(|h|^2_\infty+|g|^2_\infty)|\Omega|$ and $|\Omega|$ stands for the Lebesgue measure of $\Omega$. Next, testing Eq. $(iv)$ in Definition 2.1 by $-\Delta v$ and Eq. $(v)$ in Definition 2.1 by $-\Delta w$, and summing up the resulting equalities
we see that
\begin{equation} \frac{d}{dt}|\nabla v|^2_H+\frac{d}{dt}|\nabla w|^2_H+|\Delta v|^2_H +|\Delta w|^2_H \leq C_1   \label{3.5}\end{equation}

a.e. on $(0,T)$. From  Definition 2.1 $(iii)(a), (b)$ it follows that
\begin{equation} \sigma'=\left\{
\begin{array}{cl}
 F(\sigma,v,w)u +av' & \hspace{0.3cm} \mbox{if}  \hspace{0.3cm} f_*(v,w)< \sigma < f^*(v,w) ,\\
 f_{*v}'(v,w) v'+f_{*w}'(v,w) w' & \hspace{0.3cm} \mbox{if}  \hspace{0.3cm}  \sigma = f_*(v,w),\\
 {f_v^*}'(v,w) v' +{f_w^*}'(v,w) w' & \hspace{0.3cm} \mbox{if}  \hspace{0.3cm} \sigma =f^*(v,w) .\\
\end{array}
\right. \label{3.6}
\end{equation}

Multiplying the first line of $(\ref{3.6})$ by $\sigma'$ with the help of Young's inequality we deduce that
\begin{equation} |\sigma'|^2_H  \leq C_2 \left(1+|v'|^2_H \right)   \label{3.7}\end{equation}

a.e. on $(0,T)$ with $C_2=2\max\{m^2|F|^2_\infty|\Omega|,a^2\}$. Hence, from $(\ref{3.6})$ and $(\ref{3.7})$ we see that always
\begin{equation} |\sigma'|^2_H  \leq C_3 \left(1+|v'|^2_H+ |w'|^2_H \right)   \label{3.8}\end{equation}

a.e. on $(0,T)$, where $C_3=C_2+\max\{|f_{*v}'|_\infty, |f_{*w}'|_\infty, |{f_v^*}'|_\infty, |{f_w^*}'|_\infty\}$.

\hspace{0.6cm}  Calculating $(\ref{3.4})$+$(\ref{3.5})$+$\frac{1}{2C_3}\times$ $(\ref{3.8})$ we obtain
\begin{align*}\nonumber |v'|^2_H +|w'|^2_H +\frac{1}{C_3}|\sigma'|^2_H&+2|\Delta v|^2_H +2|\Delta w|^2_H \\ &+4\frac{d}{dt}\left\{|\nabla v|^2_H+|\nabla w|^2_H\right\}\leq 4C_1+1.  \label{}\end{align*}

Integrating now this inequality from $0$
 to $T$ we obtain the uniform with respect to $u$ estimates $(\ref{3.3})$.
 \epf

\hspace{0.6cm} Let  $\mathcal{L}:S_m\to C([0,T];H\times H\times H)$ be the operator  which
with each $u\in S_m$ associates the unique solution
\begin{equation}
\{\sigma(u),v(u),w(u)\}=\mathcal{L}(u). \label{3.9}
\end{equation}

of system  $(\ref{1.1})$--$(\ref{1.5})$. Then, we have the following result.

\begin{theorem} \label{th3-2}
The solution operator $\mathcal{L}:S_m\to C([0,T];H\times H\times H)$
is weak-strong continuous.
\end{theorem}

\bpf{Proof} The set $S_m$ endowed with the weak topology of
the space $L^2(0,T;H)$ is metrizable. Hence, it is enough to establish the
sequential continuity of the operator $\mathcal{L}$. To this aim, take  an arbitrary sequence  $u_n$,
$n\geq 1$, from $S_m$ which weakly converges  to
some $u\in S_m$. Let $\{\sigma(u_n),v(u_n),w(u_n)\}$, $n\geq 1$,
be the sequences of solutions of system  $(\ref{1.1})$--$(\ref{1.5})$
corresponding to the controls $u_n$, $n\geq 1$. By the weak and weak-star compactness
results, the uniform estimates (\ref{3.2}), (\ref{3.3}) imply that there exists a subsequence $\{\sigma(u_{n_k}):=\sigma_k,v(u_{n_k}):=v_k,w(u_{n_k}):=w_k\}$, $k\geq 1$, of
the sequence $\{\sigma(u_n),v(u_n),w(u_n)\}$, $n\geq 1$, and some elements $\sigma\in W^{1,2}(0,T;H)$, $v,w\in W^{1,2}(0,T;H)\cap
L^\infty(0,T;V)\cap L^2(0,T;H^2(\Omega))$ such
that
\begin{align}
& \left. \begin{array}{lc} &\*\hspace{-1.7cm} {v}_k \to v \;\; \mbox{and} \;\; {w}_k \to w \quad
  \vspace{0.1cm}\mbox{weakly-star in } L^{\infty}(0,T; V)\\
& \*\hspace{3cm} \mbox{and weakly in } W^{1,2}(0,T; H) \cap
L^2(0,T;H^2(\Omega))\;\;
  \vspace{0.1cm}\\
& \*\hspace{3.5cm}\mbox{and, thus, strongly in } C([0,T];H),
\end{array}
\right. \label{3.10}\\
&\*\hspace{1.2cm}{\sigma}_k \to \sigma \quad \mbox{weakly in }
W^{1,2}(0,T; H) .\label{3.11}
\end{align}

\hspace{0.6cm} Next, in order to justify the passage to the limit in the nonlinear right-hand sides of system  (\ref{1.1})--(\ref{1.3}) we show that along with
$(\ref{3.11})$ we have the following convergence
\begin{equation}{\sigma}_k \to \sigma \quad \mbox{strongly in } C([0,T];H).
\label{3.12}
\end{equation}

To this end, define the function
\begin{equation}M(\sigma,v,w):=\sigma-[\sigma-f^*(v,w)]^++[f_*(v,w)-\sigma]^+,
\label{3.13}
\end{equation}

$\sigma,v,w\in \mathbb{R}$, where $[\cdot]^+$ is the positive part of a function. Take $i,j\ge 1$ such that $i\neq j$. Then,
$$f_*(v_j,w_j)\leq M(\sigma_i,v_j,w_j).$$

In fact, if $\sigma_i<f_*(v_j,w_j)$ ($\sigma_i>f^*(v_j,w_j)$), then $M(\sigma_i,v_j,w_j)=f_*(v_j,w_j)$ $(f^*(v_j,w_j))\geq f_*(v_j,w_j).$
When $f_*(v_j,w_j)\leq\sigma_i\leq f^*(v_j,w_j)$, then $M(\sigma_i,v_j,w_j)$ $=\sigma_i\geq f_*(v_j,w_j)$ by the assumption. Similarly, we have
$$ M(\sigma_i,v_j,w_j)\leq f^*(v_j,w_j),$$

so that $M(\sigma_i,v_j,w_j)\in \mathcal{K}(v_j,w_j)$. Hence, $I_{\mathcal{K}(v_j,w_j)}(M(\sigma_i,v_j,w_j))=0$ and from the definition of the subdifferential $\partial I_{\mathcal{K}(v,w)}$
we deduce that the zero element of the space $H$
\begin{equation}\Theta_H\in\partial I_{\mathcal{K}(v_j,w_j)}(M(\sigma_i,v_j,w_j)).
\label{3.14}
\end{equation}

The monotonicity of the operator $\partial I_{\mathcal{K}(v,w)}$, Eq. $(iii)$ of Definition 2.1, and $(\ref{3.14})$ then imply that
\begin{equation}\left<\sigma_j-M(\sigma_i,v_j,w_j),F(\sigma_j,v_j,w_j)-\sigma_{j}'+av_j'-\Theta_H\right>_H\geq 0
\label{3.15}
\end{equation}

 a.e. on $[0,T]$. Furthermore, from  $(\ref{3.13})$, Definition 2.1 $(iii)(a)$, and the Lipschitz continuity of the functions $f_*$ and $f^*$ it follows that
\begin{align}\nonumber|\sigma_i-M(\sigma_i,v_j,w_j)|&=|[f_*(v_j,w_j)-\sigma_i]^+-[\sigma_i-f^*(v_j,w_j)]^+| \\ &\leq L_0(|v_j-v_i|+|w_j-w_i|)
\label{3.16}
\end{align}

a.e. on $Q(T)$, where $L_0$ is a common Lipschitz constant of $f_*$ and $f^*$. From  $(\ref{3.15})$ and  $(\ref{3.16})$ we conclude that
\begin{align*}\nonumber \langle\sigma_j&-\sigma_i,\sigma_{j}'-F(\sigma_j,v_j,w_j)-a v_j'\rangle_H \\ &\leq L_0\langle|v_j-v_i|+|w_j-w_i|,|\sigma_{j}'-F(\sigma_j,v_j,w_j)-a v_j'|\rangle_H
\label{}
\end{align*}

a.e. on $[0,T]$. Interchanging the roles of the indices $i$ and $j$ we also have
\begin{align*}\nonumber \langle\sigma_i&-\sigma_j,\sigma_{i}'-F(\sigma_i,v_i,w_i)-a v_i'\rangle_H \\ &\leq L_0\langle|v_j-v_i|+|w_j-w_i|,|\sigma_{i}'-F(\sigma_i,v_i,w_i)-a v_i'|\rangle_H
\label{}
\end{align*}

a.e. on $[0,T]$. Summing the last two inequalities up from H\"{o}lder's inequality we obtain
\begin{align*}\nonumber \langle\sigma_j&-\sigma_i,\sigma_j'-\sigma_i'\rangle_H \leq \langle\sigma_j-\sigma_i,F(\sigma_j,v_j,w_j)-F(\sigma_i,v_i,w_i)\rangle_H \\ &+ a\langle\sigma_j-\sigma_i, v_j'- v_i'\rangle_H  \\ &+6L_0\left(|\sigma_{j}'|_H+|\sigma_{i}'|_H+|a|\left(|v_{j}'|_H+|v_{i}'|_H\right)+R_1\right)\left(|v_j-v_i|_H+|w_j-w_i|_H\right),
\label{}
\end{align*}

a.e. on $[0,T]$, where $R_1=2|F|_\infty|\Omega|^\frac{1}{2}$. The application of Young's inequality further gives
\begin{align*}\nonumber \frac{d}{dt}|\sigma_j&-\sigma_i|^2_H \leq R_2\left(|\sigma_j-\sigma_i|_H^2+|v_j-v_i|_H^2+|w_j-w_i|^2_H\right) \\ &+ 2|a|\langle\sigma_j-\sigma_i, v_j'- v_i'\rangle_H   \\ &+12L_0\left(|\sigma_{j}'|_H+|\sigma_{i}'|_H+|a|\left(|v_{j}'|_H+|v_{i}'|_H\right)+R_1\right)\left(|v_j-v_i|_H+|w_j-w_i|_H\right),
\label{}
\end{align*}

a.e. on $[0,T]$, where $R_2=3+L_0^{2}$. Integrating this inequality  from $0$ to $t\in [0,T]$ we infer that
\begin{align}\nonumber |\sigma_{j}-\sigma_{i}|_H^2(t) &\leq    R_2\int_0^t |\sigma_j-\sigma_i|^2_H(\tau)\, d\tau \\ \nonumber &+R_2T\left(|v_j-v_i|^2_{C([0,T];H)}+|w_j-w_i|^2_{C([0,T];H)}\right) \\ \nonumber & +2|a| \int_0^t\int_\Omega(\sigma_j-\sigma_i)(\tau)(v_j'-v_i')(\tau) \,dx\, d\tau\\ & + R_3\left(|v_j-v_i|_{C([0,T];H)}+|w_j-w_i|_{C([0,T];H)}\right),  \label{3.17}\end{align}

$t\in [0,T]$, where $R_3=12 L_0 (2R_0^2(1+|a|)+R_1T)$. By applying Fubini's theorem and then integrating by parts, the second integral on the right-hand side of $(\ref{3.17})$ can be rewritten  and evaluated as follows
\begin{align}\nonumber 2|a|& \int_\Omega(v_j-v_i)(t) (\sigma_j-\sigma_i)(t)\,dx-
2|a| \int_\Omega\int_0^t(v_j-v_i)(\tau) (\sigma_j'-\sigma_i')(\tau)\,d\tau\, dx \\ \nonumber &\leq
2|a|  |v_j-v_i|_H(t) |\sigma_j-\sigma_i|_H(t)+
2|a| \int_0^t|v_j-v_i|_H(\tau) |\sigma_j'-\sigma_i'|_H(\tau)\,d\tau
\\ &\leq
R_4  |v_j-v_i|_{C([0,T];H)}, \label{3.18}\end{align}

where $R_4=4|a|R_0(|\Omega|^\frac{1}{2}+1)$. Therefore, applying Gronwall's inequality to $(\ref{3.17})$ in view of $(\ref{3.18})$
 and the convergences $(\ref{3.10})$ we conclude that $\sigma_k$, $k\geq 1$, is a Cauchy sequence in the space $C([0,T];H)$. Hence, according
to $(\ref{3.11})$ we obtain the convergence $(\ref{3.12})$.

\hspace{0.6cm} Now, the Lipschitz continuity of $F,g,h$ and the convergences $(\ref{3.10})$, $(\ref{3.12})$ allow us to conclude that
 \begin{align}
& \left. \begin{array}{cl}
& F(\sigma_k,v_k,w_k)\to F(\sigma,v,w), \quad h(\sigma_k,v_k,w_k)\to h(\sigma,v,w),\;\;
  \vspace{0.1cm}\\
& g(\sigma_k,v_k,w_k)\to g(\sigma,v,w) \quad \mbox{in } \, C([0,T];H)
\end{array}
\right. \label{3.19}
\end{align}
We thus also have
\begin{equation} h(\sigma_k,v_k,w_k)u_k\to h(\sigma,v,w)u \quad \mbox{weakly in } L^{2}(0,T; H)
\label{3.20}
\end{equation}

and
\begin{equation} F(\sigma_k,v_k,w_k)+av_k'-\sigma_k'\to F(\sigma,v,w)+av'-\sigma' \quad \mbox{weakly in } L^{2}(0,T; H).
\label{3.21}
\end{equation}

\hspace{0.6cm} Given the convergences $(\ref{3.10})$, $(\ref{3.11})$, and $(\ref{3.19})$--$(\ref{3.21})$ to finish the proof and show that the triple $\{\sigma,v,w\}$ is a solution to $(\ref{1.1})$--$(\ref{1.5})$ with $u\in S_m$, i.e.
\begin{equation*} \{\sigma,v,w\}=\mathcal{L}(u)= \{\sigma(u),v(u),w(u)\},
\label{}
\end{equation*}

it remains to show that
\begin{equation}  F(\sigma,v,w)u+av'-\sigma' \in \partial I_{\mathcal{K}(v,w)}(\sigma)
\label{3.22}
\end{equation}

a.e. on $[0,T]$. To this end, take an arbitrary $z\in L^{2}([0,T]; H)$ such that $z\in \mathcal{K}(v,w)$ a.e. on $[0,T]$ and for every $k\ge 1$ define the function
\begin{equation*}z_k:=z-[z-f^*(v_k,w_k)]^++[f_*(v_k,w_k)-z]^+.
\label{}
\end{equation*}

Then, as above we see that $z_k\in \mathcal{K}(v_k,w_k)$ a.e. on $[0,T]$ and from the definition of $\mathcal{K}(v,w)$ and $(\ref{3.10})$ it also follows that
\begin{equation} z_k\to z \quad \mbox{in }  L^{2}([0,T]; H).
\label{3.23}
\end{equation}

 Consequently, the definition and monotonicity of the operator $\partial
I_{\mathcal{K}(v_k,w_k)}$ imply in view of Eq. $(iii)$ of Definition 2.1 that
$$\langle F(\sigma_k,v_k,w_k)u_k+av_k'-\sigma_k',z_k-\sigma_k\rangle_H\leq 0, \quad k\geq 1, $$

a.e. on $[0,T].$ Passing in this inequality to the limit as $k\to\infty$ we see from (\ref{3.12}), (\ref{3.21}), (\ref{3.23}) that
$$\langle F(\sigma,v,w)u+av'-\sigma',z-\sigma\rangle_H\leq 0 $$

a.e. on $[0,T]$ for any $z\in L^2(0,T; H)$, $z\in \mathcal{K}(v,w)$, and thus
(\ref{3.22}) follows.

\hspace{0.6cm}  Therefore, $\{\sigma,v,w\}=\mathcal{L}(u)$ and from the
uniqueness of a solution to $(\ref{1.1})$--$(\ref{1.5})$ coupled with
the convergences (\ref{3.10}), (\ref{3.12}) it follows that $\mathcal{L}(u_n)\to
\mathcal{L}(u)$ in $C([0,T], H\times H\times H)$ hence proving the assertion of the theorem.
\epf

\hspace{0.6cm}  The next theorem provides further continuity properties of the operator $\mathcal{T}$ which are instrumental for the proof of both existence and relaxation for our control problem in the next section. To prove this theorem we will require the following lemma.

\begin{lemma}\textup{(\hspace{-0.001cm}\cite[Lemma 3.8]{AikiMinchev})}  \label{parabolic}
Let $\mu_0>0$ and $\theta$ be a solution of the initial boundary value problem
\begin{align*}&\theta'-\mu_0\Delta \theta=f \quad \text{in } Q,\\
&\frac{\partial\theta}{\partial n}=0  \quad \text{on } (0,T)\times\partial \Omega, \;\; \theta(0)=\theta_0,\end{align*}

where $f$ and $\theta_0$ are given functions. If $f\in L^r(0,T;L^q(\Omega))$ with $\frac{1}{r}+\frac{N}{2q}<1$ for $q,r\geq 1$ and
$\theta_0\in L^\infty(\Omega)$, then there exists a positive constant $C_*$ depending on $\Omega$, $\mu_0$, $q$, $r$, and $N$ only such that
$$|\theta|_{L^\infty(0,t;L^\infty(\Omega))}\leq C_*\left(|f|_{L^r(0,t;L^q(\Omega))}+|\theta_0|_{L^\infty(\Omega)}\right) \quad \text{for } 0\leq t\leq T.$$
\end{lemma}

\begin{theorem} \label{Lipschitz}
Let $u_i \in S_m$  and $\{\sigma_i,v_i, w_i\} = \mathcal{L}(u_i)$, $i = 1, 2$. Then,
\begin{align}
 |\sigma_1(t)-\sigma_2(t)|_H^2+ |v_1(t) - v_2(t)|_H^2  &+ |w_1(t)-w_2(t)|_H^2  \nonumber \\
&\leq  \  C_m   \int_0^t |u_1(\tau)-u_2(\tau)|_H^2\,d\tau,  \label{3.24}
\end{align}

$t\in [0,T]$, for a constant $C_m>0$ which depends on $m$ only.
\end{theorem}

\bpf{Proof} Denote $\sigma := \sigma_1 -\sigma_2$, $v := v_1 -v_2$, $w := w_1 - w_2$, and $u := u_1 - u_2$. Taking the difference of Eqs. $(iv)$ in Definition 2.1 corresponding to $u_1$ and $u_2$, and testing the result by $v$ we obtain
\begin{equation*} \frac{1}{2}\frac{d}{dt}|v|^2_H+|\nabla v|^2_H \leq \left|h(\sigma_1,v_1,w_1) u_1\pm h(\sigma_1,v_1,w_1) u_2 - h(\sigma_2,v_2,w_2)u_2\right|_H|v|_H. \label{}\end{equation*}

Since the second term on the left-hand side of this inequality is always nonnegative, invoking Young's inequality and the Lipschitz continuity of $h$ we have
\begin{equation*} \frac{d}{dt}|v|^2_H \leq |h|^2_\infty |u|_H^2+ 2|v|_H^2+m^2L^2\left(|\sigma|_H^2+|v|_H^2+|w|_H^2\right). \label{}\end{equation*}

The integration from $0$ to $t$ further yields
\begin{equation*}  |v(t)|_H^2\leq C_0 \int_0^t\left(|\sigma(\tau)|_H^2+|v(\tau)|_H^2+|w(\tau)|_H^2+|u(\tau)|_H^2\right)d\tau,
 \end{equation*}

$t\in[0,T]$, where $C_0=|h|^2_\infty+m^2L^2+2$. Making use of Eq. $(v)$ in Definition 2.1 a similar inequality can be obtained for $|w(t)|_H^2$  so that we have
\begin{equation} |v(t)|_H^2+|w(t)|_H^2\leq C_1 \int_0^t\left(|\sigma(\tau)|_H^2+|v(\tau)|_H^2+|w(\tau)|_H^2+|u(\tau)|_H^2\right)d\tau,
\label{3.25} \end{equation}

$t\in[0,T]$, where $C_1=2C_0$.

\hspace{0.6cm}  Multiplying now the result of the substraction of Eq. $(iv)$ in Definition 2.1 corresponding to $u_2$ from that for $u_1$ by $v'$ we see that
\begin{equation*}|v'|^2_H + \frac{d}{dt}|\nabla v|^2_H\leq |h|_\infty |u|_H|v'|_H+m|h(\sigma_1,v_1,w_1) - h(\sigma_2,v_2,w_2)|_H|v'|_H. \label{}\end{equation*}

Applying Young's inequality to the last inequality, using the Lipschitz continuity of $h$, and integrating over $(0,t)$, $t\in [0,T]$, we obtain
\begin{equation}\int_0^t|v'(\tau)|_H^2 d\tau \leq C_2 \int_0^t\left(|\sigma(\tau)|_H^2+|v(\tau)|_H^2+|w(\tau)|_H^2+|u(\tau)|_H^2\right)d\tau,
\label{3.26} \end{equation}

$t\in[0,T]$, where $C_2=2(|h|_\infty^2+m^2L^2)$.

\hspace{0.6cm} The application of Gronwall's inequality to (\ref{3.25}) leads to
\begin{equation}  |v(t)|_H^2+|w(t)|_H^2\leq C_3 \int_0^t \left(|\sigma(\tau)|_H^2+|u(\tau)|_H^2\right)d\tau,
\label{3.27} \end{equation}

$t\in[0,T],$ for $C_3=\exp\{C_1T\}$, and hence from (\ref{3.26}) we infer that
\begin{equation}\int_0^t|v'(\tau)|_H^2 d\tau \leq C_4 \int_0^t\left(|\sigma(\tau)|_H^2+|u(\tau)|_H^2\right)d\tau,
\label{3.28} \end{equation}

$t\in[0,T],$ where $C_4=C_2(1+C_3T)$. From (\ref{3.27}) we also deduce that
\begin{equation}  |v|_{L^\infty(0,t;H)}^2+|w|_{L^\infty(0,t;H)}^2\leq C_3 \int_0^t \left(|\sigma(\tau)|_H^2+|u(\tau)|_H^2\right)d\tau,
\label{3.29} \end{equation}

$t\in[0,T]$.

\hspace{0.6cm}   Now for $s\in (0,T]$ define
\begin{align}\nonumber l(s): = \max\{ &|f_*(v_1,w_1)-f_*(v_2,w_2)|_{L^{\infty}(0,s; L^{\infty}(\Omega))}, \\ &|f^*(v_1,w_1)-f^*(v_2,w_2)|_{L^{\infty}(0,s; L^{\infty}(\Omega))}\} \label{3.30} \end{align}
and
\begin{align*} \tilde{\sigma}_1: = \sigma_1 - [\sigma - l(s)]^+, \quad
\tilde{\sigma}_2: = \sigma_2 + [\sigma - l(s)]^+
\end{align*}

a.e. on $(0,s)\times \Omega$. Then, it is easily verified that the functions $\tilde{\sigma}_1$ and $\tilde{\sigma}_2$ can be taken as $z$ in Definition 2.1 (iii)(b) and thus we have
\begin{align*}
\big(\sigma'_1(t),[\sigma(t) - l(s)]^+\big)_H\leq \left(F(\sigma_1(t),v_1(t),w_1(t))+av_1',[\sigma(t) - l(s)]^+\right)_H
\end{align*}
and
\begin{align*}
-\big(\sigma'_2(t),[\sigma(t) - l(s)]^+\big)_H \leq -\left(F(\sigma_2(t),v_2(t),w_2(t))+av_2',[\sigma(t) - l(s)]^+\right)_H
\end{align*}

for a.e. $t\in [0,s]$.
Adding the last two inequalities we get
\begin{align*}
\big(\sigma'(t)&,[\sigma(t) - l(s)]^+\big)_H  \\
&\leq \left(F(\sigma_1(t),v_1(t),w_1(t)) - F(\sigma_2(t),v_2(t),w_2(t)),[\sigma(t) - l(s)]^+\right)_H \\
&+ a\left(v'(t),[\sigma(t) - l(s)]^+\right)_H
\end{align*}

for a.e. $t\in [0,s]$. The Lipschitz continuity of $F$ and Young's inequality further imply that
\begin{align*}
\frac{d}{dt}  \big|[\sigma(t)& - l(s)]^+\big|^2_H \\&\leq
 L^2  \left(|\sigma(t)|_H^2+|v(t)|_H^2+|w(t)|_H^2\right)+a^2|v'(t)|_H^2+ 2\left|[\sigma(t) - l(s)]^+\right|_H^2
 \end{align*}

for a.e. $t\in [0,s]$. From Gronwall's inequality it then follows that
\begin{align}
  \big|[\sigma(t) - l(s)]^+\big|^2_H
   \leq
C_5 \int_0^t\left(|\sigma(\tau)|_H^2+|v(\tau)|_H^2+|w(\tau)|_H^2+|v'(\tau)|_H^2\right)d\tau,
\label{3.31} \end{align}
$t\in [0,s]$, where $C_5=(L^2+a^2)\exp\{2T\}$.  To estimate the right-hand side of (\ref{3.31}) we use (\ref{3.27}), (\ref{3.28}) and obtain
\begin{align*} \left|[\sigma(t)-l(s)]^+\right|_H^2 \leq  C_6 \int_0^t \left(|\sigma(\tau)|_{H}^2+|u(\tau)|_{H}^2\right) d\tau
\end{align*}

$t\in[0,s]$,  where $C_6=C_5[(1+C_3T)+C_4]$. We can obtain a similar estimate  for $[-\sigma(t) - l(s)]^+$, and thus we have
\begin{align}
  \big|[\sigma(t) - l(s)]^+\big|^2_H+\big|[-\sigma(t) - l(s)]^+\big|^2_H\leq
  C_7 \int_0^t \left(|\sigma(\tau)|_{H}^2 + |u(\tau)|_{H}^2\right) d\tau,
\label{3.32}  \end{align}

$t\in [0,s]$, for some constant $C_7>0$.

\hspace{0.6cm} From Lemma 3.1 and the Lipschitz continuity of $h$ we infer that
\begin{align*} |v|_{L^\infty(0,t;L^\infty(\Omega))} &\leq  C_*|h(\sigma_1, v_1, w_1)-h(\sigma_2, v_2, w_2)|_{L^8(0,t;H)}
\\
& \leq  C_*L\left(|\sigma|_{L^8(0,t;H)}+|v|_{L^8(0,t;H)}+|w|_{L^8(0,t;H)}  \right),
\label{} \end{align*}

$t\in[0,T]$. Similarly, we see that
\begin{align*} |w|_{L^\infty(0,t;L^\infty(\Omega))}\leq  C_*L\left(|\sigma|_{L^8(0,t;H)}+|v|_{L^8(0,t;H)}+|w|_{L^8(0,t;H)}  \right),
\label{} \end{align*}

\hspace{0.6cm}  From these two inequalities, $(\ref{3.30})$ and $(\ref{3.29})$ we deduce that
\begin{align}\nonumber l^2(s)&\leq 4L^2\left(|v|_{L^\infty(0,s;L^{\infty}(\Omega))}^2+|w|_{L^\infty(0,s;L^{\infty}(\Omega))}^2\right) \\ &\leq  C_8\int_0^s \left(|\sigma(\tau)|_{H}^2 + |u(\tau)|_{H}^2\right) d\tau + C_9 |\sigma|^2_{L^8(0,s;H)},
\label{3.33} \end{align}

$s\in[0,T]$,  where $C_8=C_9 C_3 T^{\frac{1}{4}}$, $C_9=24L^4 C_*^2$. It is easy to see that
\begin{align}\nonumber |\sigma|&\leq\big|[\sigma - l(s)]^+-[-\sigma - l(s)]^+-\sigma\big|+\big|[\sigma - l(s)]^+\big|+\big|[-\sigma - l(s)]^+\big|\\
&\leq l(s)+\big|[\sigma - l(s)]^+\big|+\big|[-\sigma - l(s)]^+\big|  \label{3.34} \end{align}

a.e. on $(0,s)\times \Omega$. Therefore, from $(\ref{3.32})$--$(\ref{3.34})$ we conclude that
\begin{align*}
  |\sigma|_{L^\infty(0,t;H)}^2\leq
  C_{10} \int_0^t \left(|\sigma(\tau)|_{H}^2 + |u(\tau)|_{H}^2\right) d\tau+ C_{11}  |\sigma|_{L^8(0,t;H)}^2
\label{}  \end{align*}

$t\in[0,T]$, where $C_{10}=3(C_7+C_8|\Omega|)$, $C_{11}=3 C_9 |\Omega|$. Invoking H\"{o}lder's inequality from the last inequality we obtain
\begin{align*} |\sigma|_{L^\infty(0,t;H)}^2\leq  C_{12} |\sigma|_{L^8(0,t;H)}^2+ C_{10} \int_0^t |u(\tau)|_{H}^2 d\tau, \quad t\in [0,T]
\end{align*}

for $C_{12}=C_{10}T^\frac{3}{4}+C_{11}$. From this inequality we see that
\begin{align*} |\sigma(t)|_H^8\leq 6 C_{12}^4 \int_0^t|\sigma(\tau)|_{H}^8 d\tau+ 6C_{10}^4 \left(\int_0^t |u(\tau)|_{H}^2 d\tau\right)^4, \quad t\in [0,T].
\end{align*}

Gronwall's inequality then implies that
\begin{align}
  |\sigma(t)|_{H}^2\leq
  C_{13} \int_0^t  |u(\tau)|_{H}^2 d\tau,
\label{3.35}  \end{align}

$t\in[0,T]$, where $C_{13}=2 C_{10}\exp\{2 C_{12}^4 T\}$. Finally, from (\ref{3.27}) and (\ref{3.35}) we conclude that
\begin{align*} |\sigma(t)|_H^2+|v(t)|_H^2+|w(t)|_H^2  \leq C_m \int_0^t |u(\tau)|_H^2 d\tau,
\label{} \end{align*}

$t\in[0,T]$, where $C_m=C_3+C_{13}+C_3C_{13}T$. \epf

\section{Existence and relaxation for the control problem $(P)$}

In this section, we prove Theorem 2.1.  First, we establish the \emph{existence} of solutions for
the control system $(P)$. Since for a closed set $U\subset \mathbb{R}$ we evidently have $U\subset \co U$, any solution of $(P)$ is automatically a solution of $(RP)$. Theorem 3.2  implies that the image  $\mathcal{L}(S_m)$ of the set $S_m$ under the solution operator $\mathcal{L}$ is
compact in $C([0,T];H\times H\times H)$. Using the properties
$(\mathcal{U}1)$--$(\mathcal{U}3)$ of the mapping $\mathcal{U}$ it is a standard matter to show that its associated multivalued Nemytskii operator  $\Psi:\mathcal{L}(S_m)\to
L^2(0,T;H)$ defined by
\begin{equation} \*\hspace{-0.25cm}\Psi(\sigma,v,w):=\{u\in L^2(0,T;H); \; u(t)\in \mathcal{U}(t,\sigma(t),v(t),w(t)) \; \mbox{for a.e.} \; t\in[0,T]\}   \label{4.1}\end{equation}

is lower semicontinuous. Moreover, $\Psi$ has closed decomposable values. Recall that a subset of
$L^2(0,T;H)$ is called decomposable if along
with any two functions $u_1, u_2\in  L^2(0,T;H)$ it contains the
function $u_1\chi_E+u_2\chi_{[0,1]\setminus E}$ for any measurable
set $E\subset [0,T]$, where $\chi_A$ is the characteristic function of a set $A$.
Then,
\cite[Theorem 3.1]{Fryszkowski} implies that there exists a continuous selection of $\Psi$, i.e. a continuous mapping
${\psi}:\mathcal{L}(S_m)\to L^1(0,T;H)$ such that
\begin{equation} {\psi}(\sigma,v,w)\in \Psi(\sigma,v,w), \; (\sigma,v,w)\in \mathcal{L}(S_m) .\label{4.2}\end{equation}

By virtue of  $(\mathcal{U}3)$ we see that, in fact, ${\psi}$ is continuous from $\mathcal{L}(S_m)$ to $L^2(0,T;H)$
and ${\psi}(\sigma,v,w) \in S_m$, $(\sigma,v,w)\in \mathcal{L}(S_m)$.

\hspace{0.6cm} Next, consider  the superposition ${\psi}\circ \mathcal{L}$ of $\mathcal{L}$ and
${\psi}$. Theorem 3.2 implies that
${\psi}\circ \mathcal{L}:S_m\to S_m$ is weak-weak continuous.
The fact that the set $S_m$ is evidently convex and compact in the weak topology of
the space $L^2(0,T;H)$ allows us to infer, invoking the Schauder fixed point theorem,  that there exists a fixed point $u_*\in S_m$ of the
operator ${\psi}\circ \mathcal{L}$:
\begin{equation} u_*={\psi}\circ \mathcal{L}(u_*)={\psi}(\mathcal{L}(u_*)) .\label{4.3}\end{equation}

Letting $(\sigma_*,v_*,w_*):=\mathcal{T}(u_*)$, from (\ref{4.1})--(\ref{4.3})
we finally conclude  that $(\sigma_*,v_*,w_*,u_*)$ is a solution to Problem $(P)$.

\hspace{0.6cm} Now, to prove the \emph{relaxation}, take an arbitrary solution
$(\sigma_*,v_*,w_*,u_*)$  to the convexified problem $(RP)$. In particular, we have $u_*(t)\in \overline{\co}\,{\mathcal{U}}(t,\sigma_*(t),v_*(t),w_*(t))$, $t\in [0,T]$. In view of the properties $(\mathcal{U}1)$--$(\mathcal{U}3)$ from \cite[Corollary 1.1]{Chuong} it follows that for any $n\geq 1$
there exists a measurable function ${\gamma}_n(t)\in {\mathcal{U}}(t,\sigma_*(t),v_*(t),w_*(t))$, $t\in [0,T]$, such that
\begin{equation}
\sup\limits_{0\leq s\leq t\leq T} \left|\int_s^{t}
\left({u}_*(\tau)-{\gamma}_n(\tau)\right) \, d\tau\right|_H \leq
\frac{1}{n} .\label{4.4}
\end{equation}

From $(\mathcal{U}3)$ we see that for any $(\sigma,v,w)\in H\times H\times H$ and a.e. $t\in [0,T]$ there exists ${\gamma}\in {\mathcal{U}}(t,\sigma,v,w)$ such that
\begin{align*}|\gamma_n(t)-\gamma|_H^2&< \frac{2}{n^2}+2 d^{\,2}_H(\gamma_n(t),{\mathcal{U}}(t,\sigma,v,w))\\
&< \frac{2}{n^2}+12k^2(t)\left(|\sigma_*(t)-\sigma|_H^2+|v_*(t)-v|_H^2+|w_*(t)-w|_H^2\right). \label{}\end{align*}

Making use of this inequality we now construct the following multivalued mapping
\begin{align}\nonumber \mathcal{U}_n(t,\sigma,v,w)&:=\{ \gamma\in \mathcal{U}(t,\sigma,v,w); \; |\gamma_n(t)-\gamma|_H^2 \\
&\leq \frac{2}{n^2}+12k^2(t)\left(|\sigma_*(t)-\sigma|_H^2+|v_*(t)-v|_H^2+|w_*(t)-w|_H^2\right)\}, \label{4.5}\end{align}

and define its associated multivalued  Nemytskii  operator $\Psi_n:\mathcal{L}(S_m)\to
L^2(0,T;H)$ similarly as in (\ref{4.1}):
\begin{equation} \*\hspace{-0.515cm}\Psi_n(\sigma,v,w):=\{u\in L^2(0,T;H); \; u(t)\in {\mathcal{U}_n(t,\sigma(t),v(t),w(t))} \; \mbox{for a.e.} \; t\in[0,T]\}.   \label{4.6}\end{equation}

As above, invoking \cite[Theorem 3.1]{Fryszkowski} we find a continuous mapping
${\psi}_n:\mathcal{L}(S_m)\to S_m$ such that
\begin{equation} {\psi}_n(\sigma,v,w)\in \Psi_n(\sigma,v,w), \; (\sigma,v,w)\in \mathcal{L}(S_m) \label{4.7}\end{equation}

and invoking the Schauder fixed point theorem we find a fixed point $u_n\in S_m$ of the
superposition ${\psi}_n\circ \mathcal{L}$:
\begin{equation} u_n={\psi}_n(\mathcal{L}(u_n)) .\label{4.8}\end{equation}

Letting $(\sigma_n,v_n,w_n):=\mathcal{L}(u_n)$, from (\ref{4.5})--(\ref{4.8})
we conclude that $(\sigma_n,v_n,w_n,u_n)$, $n\geq 1$, is a solution to problem $(P)$ and, in addition, we have
\begin{align}\nonumber
|&{\gamma}_n(t)-u_n(t)|_H^2 \\ &\leq \frac{2}{n^2}
+12k^2(t)\left(|\sigma_*(t)-\sigma_n(t)|_H^2+|v_*(t)-v_n(t)|_H^2+|w_*(t)-w_n(t)|_H^2\right).
\label{4.9}
\end{align}

From the fact that on the set $S_m$ the weak topology of the space $L^2(0,T;H)$ coincides with the topology generated by the ``weak norm'' given by the supremum on the left-hand side of $(\ref{4.4})$, we infer, in view of $(\ref{4.4})$, that
\begin{equation}
\gamma_n\to u_* \quad \text{weakly in} \; \;L^2(0,T;H).
 \label{4.10}\end{equation}

Then, from Theorem 3.2 we obtain
\begin{align}(\sigma(\gamma_n),v(\gamma_n),w(\gamma_n)) \to (\sigma_*,v_*,w_*) \quad \mbox{strongly in }
C([0,T];H^3) . \label{4.11}
\end{align}

Combining Theorem 3.3 with $(\ref{4.9})$ we have
\begin{align*}
|\sigma_n(t)&-\sigma(\gamma_n)(t)|_H^2 +|v_n(t)-v(\gamma_n)(t)|_H^2 +|w_n(t)-w(\gamma_n)(t)|_H^2\\ &\leq
C_m\int_0^t
|u_n(\tau)-\gamma_n(\tau)|_H^2 d\tau \\
 &\leq C_m\int_0^t
\bigg(\frac{2}{n^2}+24k^2(\tau)\big(|\sigma_n(\tau)-\sigma(\gamma_n)(\tau)|_H^2
+|v_n(\tau)-v(\gamma_n)(\tau)|_H^2 \\
&\*\hspace{1.9cm}+ |w_n(\tau)-w(\gamma_n)(\tau)|_H^2\big)
+24k^2(\tau)\big(|\sigma(\gamma_n)(\tau)-\sigma_*(\tau)|_H^2 \\
&\*\hspace{1.9cm}+|v(\gamma_n)(\tau)-v_*(\tau)|_H^2 +
|w(\gamma_n)(\tau)-w_*(\tau)|_H^2\big)\bigg) d\tau.\end{align*}

This inequality, (\ref{4.11}), and Gronwall's  inequality further yield
\begin{equation*}
(\sigma_n,v_n,w_n) \to (\sigma_*,v_*,w_*) \quad \mbox{strongly in } C([0,T];H\times H\times H),
\end{equation*}

which together with (\ref{4.9}), (\ref{4.10}) implies that
\begin{equation*}
u_n \to u_* \quad \mbox{weakly in } L^2(0,T;H).
\end{equation*}

The last two convergences finally prove the relaxation part of Theorem 2.1.

\section{Conclusion}

In this paper, we have considered a nonlinear control system of PDEs arising in population dynamics of three biological species: predator, prey and food for prey. This system is subject to a nonconvex control constraint and it takes account of the situation when the dependence of the density of the food for prey on the densities of preys and predators has a hysteretic character. By exploring and then exploiting the continuity-type properties of the control-to-state solution operator we establish the existence of solutions for our control problem and obtain some relaxation properties of it.

\vspace{0.5cm}

\textbf{Acknowledgements.}  The authors want to thank the
anonymous referees for their valuable suggestions and remarks
which helped to improve the manuscript.

\bibliographystyle{amsplain}

\end{document}